\documentclass[12pt,reqno]{amsart}
\usepackage{graphicx,subfigure}
\usepackage{hyperref}
\newcommand{\vertiii}[1]{{\left\vert\kern-0.25ex\left\vert\kern-0.25ex\left\vert #1
    \right\vert\kern-0.25ex\right\vert\kern-0.25ex\right\vert}}
\hypersetup{colorlinks=true, citecolor=blue, linkcolor=red}

\numberwithin{equation}{section} \numberwithin{figure}{section}
\numberwithin{table}{section} 
{ \theoremstyle{plain}
\newtheorem{theorem}{Theorem}[section]

\newtheorem{rem}[theorem]{Remark}

}

\begin{document}

\title[On the nonlinear Schr\"{o}dinger equation in the half-space]
{One-dimensional symmetry of positive bounded solutions to the subcubic and cubic nonlinear Schr\"{o}dinger equation in the half-space     in  dimensions $N=4,5$}



\author[Christos Sourdis]{Christos Sourdis}
\address{General Lyceum of Malia, Heraklion, Crete, Greece}
\email{sourdis@uoc.gr}

\maketitle

\begin{abstract}
We are concerned with the half-space Dirichlet problem
\[\left\{\begin{array}{ll}
  -\Delta v+v=|v|^{p-1}v & \textrm{in}\ \mathbb{R}^N_+, \\
  v=c\ \textrm{on}\ \partial\mathbb{R}^N_+, &\lim_{x_N\to \infty}v(x',x_N)=0\ \textrm{uniformly in}\ x'\in\mathbb{R}^{N-1},
\end{array}\right.
\]
where $\mathbb{R}^N_+=\{x\in \mathbb{R}^N \ : \ x_N>0\}$ for some $N\geq 2$, and $p>1$, $c>0$ are constants.
It was shown recently  by Fernandez and Weth [Math. Ann. (2021)] that there exists an explicit number $c_p\in (1,\sqrt{e})$, depending only on $p$, such that for $0<c<c_p$ there are infinitely many bounded positive solutions, whereas, for $c>c_p$ there are no bounded positive solutions. They also  posed as an interesting open question whether the
  one-dimensional solution is the unique bounded positive solution   in the case where $c = c_p$.
If $N=2,  3$, we recently showed this one-dimensional symmetry property   in [Partial Differ. Equ. Appl. (2021)] by adapting some ideas from the proof of De Giorgi's conjecture in low dimensions. Here, we first focus on the case $1<p<3$ and prove this uniqueness property in dimensions $2\leq N\leq 5$. Then, for the cubic NLS, where $p=3$, we establish this for $2\leq N \leq 4$. Our approach is completely different and relies on showing that a suitable auxiliary function, inspired by a Lyapunov-Schmidt type decomposition of the solution, is a nonnegative super-solution to a Lane-Emden-Fowler equation in $\mathbb{R}^{N-1}$, for which an optimal Liouville type result is available.
\end{abstract}

\section{Introduction}
Recently in \cite{weth}, the authors studied the half-space Dirichlet problem
\begin{equation}\label{eqEq}\left\{\begin{array}{ll}
  -\Delta v+v=|v|^{p-1}v & \textrm{in}\ \mathbb{R}^N_+, \\
 & \\
  v=c\ \textrm{on}\ \partial\mathbb{R}^N_+, &\lim_{x_N\to \infty}v(x',x_N)=0\ \textrm{uniformly in}\ x'\in\mathbb{R}^{N-1},
\end{array}\right.
\end{equation}
where $\mathbb{R}^N_+=\{x\in \mathbb{R}^N \ : \ x_N>0\}$ for some $N\geq 1$, and $p>1$, $c>0$ are constants.
We note that $V(x,t)=e^{it}v(x)$ is a standing wave solution to the focusing nonlinear Schr\"{o}dinger equation with the pure odd power nonlinearity of exponent $p$.
Let us clarify here that throughout this paper solutions will be understood in the classical sense (i.e. at least of class $C^2$, and continuous up to the boundary).

If $N=1$, then the corresponding ODE has a unique positive even solution that decays to zero at infinity, it is given explicitly by the following formula
\begin{equation}\label{eqw0}
t \rightarrow w_0(t)=c_p\left[\cosh\left(\frac{p-1}{2}t\right)
\right]^{-\frac{2}{p-1}}\ \textrm{with}\ c_p=\left(\frac{p+1}{2} \right)^{\frac{1}{p-1}}=w_0(0)=\sup_{t\in \mathbb{R}}w_0(t).
\end{equation}
Still for $N=1$, it was shown in the aforementioned reference that if $0<c<c_p$ then (\ref{eqEq}) possesses exactly two positive solutions given by
\[
t\to w_0(t+t_{c,p})\ \textrm{and}\ t\to w_0(t-t_{c,p})
\]
with
\[
t_{c,p}=\frac{2}{p-1}\ln\left(\sqrt{\frac{p+1}{2c^{p-1}}}+\sqrt{\frac{p+1}{2c^{p-1}}-1} \right);
\]
 if $c=c_p$ then $w_0$ is the unique positive solution; if $c>c_p$ then there are no positive solutions. We note in passing that the above solutions play an important role in a class of boundary layer problems (see \cite{butuzov}).

If $N\geq 2$, $p>1$, and $0<c<c_p$, using variational methods, it was shown in the same reference \cite{weth} that (\ref{eqEq}) admits at least three positive bounded solutions that are geometrically distinct in the sense that they are not translates of each other in the $x'$ direction. In particular, under the further restriction that  $p+1$ is smaller than the critical Sobolev exponent in $\mathbb{R}^N$, $N\geq 2$, then (\ref{eqEq}) admits a positive bounded solution of the form
\[
x\to w_0(x_N+t_{c,p})+\tilde{v}(x) \ \textrm{with}\ \tilde{v}\in H_0^1(\mathbb{R}^N_+)\setminus \{0\} \ \textrm{nonnegative},
\]($H_0^1(\mathbb{R}^N_+)$ denotes the usual Sobolev function space with zero trace on $\partial \mathbb{R}^N_+$).

On the other hand, if $c>c_p$, $p>1$, it was shown therein that (\ref{eqEq}) has   no bounded positive solutions. This was accomplished by means of the famous sliding method \cite{berest}.

 Still in the same reference \cite{weth}, it was posed as an interesting open question whether the function $x
\to w_0(x_N)$ is
the unique bounded positive solution to (\ref{eqEq}) in the case $c = c_p$, $p>1$.  In this regard, as we explained in \cite{sourdaras}, we note that the aforementioned sliding argument  can be applied even in this case to establish that
\begin{equation}\label{eqOrder}
  w_0(x_N)<v(x), \ x\in \mathbb{R}^N_+, \ \textrm{or}\ w_0\equiv v.
\end{equation}
If $N=2, 3$, in the aforementioned reference, we were able to exclude the first scenario in (\ref{eqOrder}) by adapting some ideas from the proof of the famous De Giorgi conjecture in the plane (see \cite{bcn,gui}). More precisely, a key observation in the proof was that the convexity of the nonlinearity  implies that, in the   first case of (\ref{eqOrder}), the difference $v-w_0$ would be a positive super-solution to the linearized problem on $w_0$ (see also \cite[Ch. 1]{dupaBook}). Actually, this property brought the problem closer in spirit to that of the one-dimensional symmetry  of bounded,  stable solutions to semilinear elliptic equations in the plane (a stronger version of De Giorgi's conjecture, see  \cite{dancer,farina}).  The fact that we were dealing with the half-space and not the full space also created some technical difficulties   in applying the approach of the aforementioned references. We point out that we were able to gain one more dimension, compared to the aforementioned works, owing to the uniform  exponential decay of solutions as $x_N\to \infty$.


In the current work, we will first restrict ourselves to the subcubic regime $1<p<3$ and prove the following.
\begin{theorem}\label{thm}
  If  $1<p<3$ and $2\leq N\leq 5$, then the only positive bounded solution of (\ref{eqEq}) with $c=c_p$  is $v(x)=w_0(x_N)$, where $c_p$ and $w_0$ are as in (\ref{eqw0}).
\end{theorem}Then, for the cubic case we can show the following.

\begin{theorem}\label{thm222}
If $p=3$ and $2\leq N \leq 4$, then the same conclusion of Theorem \ref{thm} holds.
\end{theorem}

We would like to highlight that byproducts of our proofs are  an astonishingly simple proof of our aforementioned result in
\cite{sourdaras} in  two and three dimensions (see Remark \ref{rem1} below) and  a nontrivial integral estimate that holds in all dimensions (see Remark \ref{rem2} below).

Our main observation behind the proofs is that the  auxiliary function
\[
u(x')=\int_{0}^{\infty}\left(v(x',x_N)-w_0(x_N) \right)\left(-w'_0(x_N)\right)dx_N,\ \ x'\in \mathbb{R}^{N-1},
\]
is a nonnegative super-solution to the quadratic Lane-Emden-Fowler equation $\Delta u+u^2=0$ if $1<p<3$, after a suitable re-scaling, or the cubic one if $p=3$. We point out that this property is valid for all $N\geq 2$. Then, restricting ourselves to sufficiently low dimensions, we can apply a well known but seldom utilized Liouville type result
to conclude (see Appendix \ref{App} below). We note in passing that the above auxiliary function $u$ corresponds to the projection of the difference $v-w_0$ on  $w_0'$ which is in the kernel of the linearized problem on $w_0$. Lastly, let us remark that an approach of this nature has been frequently applied in various parabolic problems for obtaining Liouville type results for ancient or eternal solutions (see for instance \cite{merle}).

Unfortunately, it seems that our method of proof does not provide useful information if $p>3$ and $N\geq 4$.

For completeness, let us mention that (\ref{eqEq}) with $c=0$ is fairly well understood. Indeed, if $p+1>2$ is subcritical in the sense of the Sobolev imbedding in $\mathbb{R}^N$, $N\geq 2$, it was shown in \cite{esteban} that
there are no nontrivial solutions in $H_0^1(\mathbb{R}^N_+)$. Furthermore, it is known that there are no positive bounded solutions (we point out that this holds without  any restriction on  $p>1$, see \cite{farinaAIMS}).
Lastly, we refer to \cite{selmi} for nonexistence results
via the Morse index.

The proofs of Theorems \ref{thm} and \ref{thm222} will be given in the following section together with some related remarks. Finally, in Appendix \ref{App} we recall a well known Liouville type theorem concerning the nonexistence of positive super-solutions to the  Lane-Emden-Fowler equation in the whole space.

\section{Proofs of the main results} In this section, we will give the proofs of Theorems \ref{thm} and \ref{thm222} along with some related remarks.
\subsection{Proof of Theorem \ref{thm}}\begin{proof}
As we have already mentioned,  relation (\ref{eqOrder}) is valid. We wish to show that the second alternative is the one which holds. To this end, let us argue by contradiction and suppose
that
\begin{equation}\label{eqordinaryb}
  w_0(x_N)<v(x',x_N), \ (x',x_N)\in \mathbb{R}^N_+.
\end{equation}

As we know from \cite{sourdaras}, the difference \begin{equation}\label{eqphi}
\varphi(x',x_N)=v(x',x_N)-w_0(x_N),\ (x',x_N) \in \overline{\mathbb{R}^N_+},
\end{equation}
 furnishes a positive super-solution to the linearization of (\ref{eqEq}) on $w_0$.
Here, we will need to make use of the full  PDE that $\varphi$ satisfies, which is nonlinear and we shall recall below. Since both $v$ and $w_0$ satisfy the same PDE in (\ref{eqEq}) and are positive, we find that
\[
\begin{array}{rcl}
  \Delta(v-w_0) & = & v-v^p-w_0+w_0^p \\
    &   &   \\
   & = & (1-pw_0^{p-1})(v-w_0)-\left(v^p-w_0^p-pw_0^{p-1}(v-w_0)\right) .
\end{array}
\]
Hence, recalling the definition of $\varphi$ from (\ref{eqphi}), and applying the mean value theorem, we arrive at
\begin{equation}\label{eqsuper}
\Delta \varphi=(1-pw_0^{p-1})\varphi-\frac{p(p-1)}{2}\xi^{p-2}\varphi^2,\ \varphi> 0 \ \textrm{in} \ \mathbb{R}^N_+;\ \varphi=0\ \textrm{on} \ \partial\mathbb{R}^N_+ ,
\end{equation}where
\begin{equation}\label{eqxi}
w_0(x_N)<\xi(x)<v(x), \ x=(x',x_N)\in \mathbb{R}^N_+.
\end{equation}

Let us now consider the positive auxiliary function
\begin{equation}\label{eqaux}
u(x')=\int_{0}^{\infty}\varphi(x',x_N)Z(x_N)dx_N,\ \ x'\in \mathbb{R}^{N-1},
\end{equation}
where $\varphi$ is as in (\ref{eqphi}), and
\begin{equation}\label{eqZ}
Z(x_N)\equiv-w'_0(x_N).
\end{equation}
For future reference, we observe that differentiation of the ODE satisfied by $w_0$  yields
\begin{equation}\label{eqReferia}
  -Z''+\left(1-pw_0^{p-1}(x_N)\right)Z=0, \ x_N>0.
\end{equation}
Furthermore, using (\ref{eqw0}), we get
\begin{equation}\label{eqZ1}
Z(0)=0,\ Z>0\ \textrm{in}\ (0,\infty)\ \textrm{and}\ Z/w_0\to1\ \textrm{as}\ x_N \to \infty.
\end{equation}

We compute that
\begin{equation}\label{eqhuge}
  \begin{array}{lll}
   \Delta_{x'}u(x') & = & \int_{0}^{\infty}\Delta_{x'}\varphi(x',x_N)Z(x_N)dx_N \\
   &  &  \\
   &  =& \int_{0}^{\infty}\left(-\varphi_{x_Nx_N}+\Delta\varphi\right)Zdx_N  \\
   &  &  \\
\textrm{using}\ (\ref{eqsuper}):   & = & \int_{0}^{\infty}\left[-\varphi_{x_Nx_N}+(1-pw_0^{p-1})\varphi-\frac{p(p-1)}{2}\xi^{p-2}\varphi^2\right]Zdx_N \\
    &   &   \\
  \textrm{integrating by parts}: & = & \int_{0}^{\infty}\varphi\left[-Z''+(1-pw_0^{p-1})Z\right]dx_N-\frac{p(p-1)}{2}\int_{0}^{\infty}\xi^{p-2}\varphi^2Zdx_N \\
  &   &   \\
  &   &  +\varphi_{x_N}(x',0)Z(0)-\varphi(x',0)Z'(0)\\
 &  & \\
\textrm{via}\ (\ref{eqsuper}), (\ref{eqReferia}), (\ref{eqZ1}): & = & -\frac{p(p-1)}{2}\int_{0}^{\infty}\xi^{p-2}\varphi^2Zdx_N.
\end{array}
\end{equation}

We first consider the case where $1<p\leq 2$. Then, by virtue of (\ref{eqxi}), we obtain
\begin{equation}\label{eqxi1}
\xi^{p-2}\geq v^{p-2}\geq \|v\|_{L^\infty(\mathbb{R}^N_+)}^{p-2}\ \textrm{in}\ \mathbb{R}^N_+.
\end{equation}
On the other hand, thanks to the Cauchy-Schwarz inequality and recalling the definition of $Z$ from (\ref{eqZ}) (keeping in mind also (\ref{eqZ1})), we have
\begin{equation}\label{eqCS}
u^2(x')=\left(\int_{0}^{\infty}\varphi Zdx_N\right)^2 \leq \int_{0}^{\infty}\varphi^2Zdx_N\int_{0}^{\infty}Zdx_N=c_p\int_{0}^{\infty}\varphi^2Zdx_N.
\end{equation}
Hence, by combining (\ref{eqhuge}), (\ref{eqxi1}) and (\ref{eqCS}), we infer that
\begin{equation}\label{eqGidarass}
-\Delta_{x'}  u\geq Cu^2,\ u> 0 \ \textrm{in}\ \mathbb{R}^{N-1},
\end{equation}
for some constant $C>0$ that depends only on $p$ and the supremum of $v$.
However, in view of Theorem \ref{thmGidaras} in Appendix \ref{App}, after a simple rescaling in the above inequality, we can arrive at a contradiction, provided that we restrict $N\geq2$ so that
\[
2\leq p_{sg}(N-1) \ \iff \ N=2\ \textrm{or}\ N=3\ \textrm{or} \ 2\leq \frac{N-1}{N-3}\ \iff \ 2\leq N\leq 5.
\]

To complete the proof, it remains to consider the case $2<p<3$.
Then, the ordering (\ref{eqxi}) gives
\begin{equation}\label{eqxi2}
\xi^{p-2}> w_0^{p-2}\ \textrm{in}\ \mathbb{R}^N_+.
\end{equation}
Now, we observe that
\[\begin{array}{rcl}
    u^2(x') & = & \left(\int_{0}^{\infty}\varphi Zdx_N\right)^2 \\
     &  &  \\
   \textrm{by the Cauchy-Schwarz inequality}:  & \leq &\int_{0}^{\infty}\frac{Z}{\xi^{p-2}}dx_N \int_{0}^{\infty}\xi^{p-2}\varphi^2Zdx_N \\
     &  &  \\
  \textrm{using}\  (\ref{eqxi2}) :    & \leq & \int_{0}^{\infty}\frac{Z}{w_0^{p-2}}dx_N \int_{0}^{\infty}\xi^{p-2}\varphi^2Zdx_N \\
   &  & \\
 \textrm{via}\ (\ref{eqZ1}) :   &\leq   & \|\frac{Z}{w_0}\|_{L^\infty(0,\infty)}\int_{0}^{\infty}w_0^{3-p}dx_N \int_{0}^{\infty}\xi^{p-2}\varphi^2Zdx_N \\
  &  & \\
 \textrm{since}\ p<3\ \textrm{and}\ w_0\ \textrm{has exponential decay} : &=  &  C' \int_{0}^{\infty}\xi^{p-2}\varphi^2Zdx_N,
  \end{array}\]
for some positive constant $C'$ that depends only on $p$. Thus, by combining (\ref{eqhuge}) and the above relation, we see that $u$ satisfies (\ref{eqGidarass}), for possibly a different positive constant, and can conclude as before.
\end{proof}
\begin{rem}\label{rem1}The approach of the current paper can be used to give a streamlined proof of our previously mentioned result in \cite{sourdaras}.
 Indeed, it follows from (\ref{eqhuge}) that $u$ is a nonnegative super-harmonic function in $\mathbb{R}^{N-1}$, $N\geq 2$, for any $p>1$. Therefore, if $N=2,3$, by the famous theorem of Hadamard and Liouville (see for instance \cite[Thm. 3.1]{farina}), we  deduce that $u$ is a constant function. In turn,  using (\ref{eqhuge}) once more, we conclude that $\varphi\equiv 0$, as desired. It is worth mentioning that an argument of this nature can be found in \cite[Lem. 6.1]{delPino} for the study of the linearized Allen-Cahn equation.
   \end{rem}
\begin{rem}\label{rem2}
We  recall from the proof of Theorem \ref{thm} that the auxiliary function $u$, as defined in (\ref{eqaux}), satisfies (\ref{eqGidarass}) for all  $N\geq 2$, provided that
$1<p<3$. Hence, by applying \cite[Lem. 2.4]{serinzou}, for any $\gamma \in (0,2)$, we infer that there exists a positive constant $M>0$ such that
\[
\int_{|x'|<R}^{}u^\gamma dx'\leq M R^{N-1-2\gamma}, \ R>0.
\]
\end{rem}\subsection{Proof of Theorem \ref{thm222}}
\begin{proof}
The proof proceeds as that of Theorem \ref{thm}, invoking a contradiction argument. However,
instead of using the mean value theorem, we can compute explicitly that the positive auxiliary function $u$, as defined in (\ref{eqaux}), satisfies
\[
\Delta_{x'}u=-\int_{0}^{\infty}(3w_0\varphi^2+\varphi^3)Zdx_N\leq -\int_{0}^{\infty}\varphi^3Zdx_N,\ x'\in \mathbb{R}^{N-1}.
\]

Next, by H\"{o}lder's inequality we get
\[
u^3(x')=\left(\int_{0}^{\infty}\varphi Zdx_N\right)^3\leq \left(\int_{0}^{\infty}Zdx_N\right)^2\int_{0}^{\infty}\varphi^3Zdx_N=c_3^2\int_{0}^{\infty}\varphi^3Zdx_N,
\]
$x'\in \mathbb{R}^{N-1}$ (recall the definition of $Z$ from (\ref{eqZ})).

Hence,   we deduce that
\[
-\Delta_{x'}u\geq \frac{1}{c_3^2}u^3,\ u>0,\ x'\in \mathbb{R}^{N-1}.
\]
On the other hand, this contradicts the Liouville type theorem in Appendix \ref{App}  since it is easy to check that
\[
  3\leq p_{sg}(N-1)\ \Longleftrightarrow \ N=2\ \textrm{or}\ N=3\ \textrm{or}\ 3\leq\frac{N-1}{N-3} \ \Longleftrightarrow \ 2\leq N\leq 4,
\]
and the proof is completed.
\end{proof}

\appendix
\section{A Liouville type theorem}\label{App}
For the reader's convenience, we state below a well known result due to \cite{gidas} (see also \cite{sirakov} and \cite[Ch. I]{qs} for  simpler proofs and extensions), which we used in the proof of Theorem \ref{thm}.

\begin{theorem}\label{thmGidaras}
Let $1 < p \leq p_{sg}(n)$, where
\[
p_{sg}(n)=\left\{\begin{array}{ll}
                   \infty & \textrm{if}\ n\leq 2, \\
                     &  \\
                   \frac{n}{n-2}  & \textrm{if}\ n> 2.
                 \end{array}
 \right.
\]
 Then, the inequality
\[-\Delta u \geq  u^p,\ x\in  \mathbb{R}^n, \]
does not possess any positive classical solution.
\end{theorem}\begin{rem}\label{remGidi}
As was remarked in \cite{qs}, the condition $p \leq p_{sg}$ in Theorem \ref{thmGidaras}  is optimal, as shown by the explicit
example $u(x) = k\left(1 + |x|^2\right)^{-1/(p-1)}$ with $n \geq  3$, $p > p_{sg}$ and $k > 0$ small enough.
\end{rem}

\end{document}